\documentclass{article}

\usepackage{amssymb,amsfonts}

\usepackage{url}

\def\cqfd{\skip10=\parfillskip\parfillskip=0pt
\enspace\hfill\symbolecqfd\par\parfillskip=\skip10\par\medskip}
\def\symbolecqfd{\rlap{$\sqcap$}$\sqcup$}

\newtheorem{theorem}{Theorem}[section]
\newtheorem{proposition}[theorem]{Proposition}
\newtheorem{lemma}[theorem]{Lemma}
\newtheorem{corollary}[theorem]{Corollary}

\newtheorem{pro-fact}[theorem]{Fact}

\newtheorem{pro-example}[theorem]{Example}
\newenvironment{example}{\begin{pro-example}\rm}{\cqfd\end{pro-example}}

\newtheorem{pro-remark}[theorem]{Remark}
\newenvironment{remark}{\begin{pro-remark}\rm}{\cqfd\end{pro-remark}}

\newenvironment{preuve}{\rm \trivlist \item[\hskip \labelsep{\bf
Proof.}]}{\cqfd\endtrivlist}

\def\cqfd{\skip10=\parfillskip\parfillskip=0pt
\enspace\hfill\symbolecqfd\par\parfillskip=\skip10\par\medskip}
\def\symbolecqfd{\rlap{$\sqcap$}$\sqcup$}
\def\proof{\begin{preuve}}
\def\eop{\end{preuve}}

\def\inter[#1]{[\![#1]\!]}
\def\inv{^{-1}}

\let\phi\varphi
\let\epsilon\varepsilon

\def\calB{\mathcal{B}}
\def\calA{\mathcal{A}}
\def\O {\mathcal{O}}

\def\leff{\le_{\sf ff}}
\def\lefg{\le_{\sf fg}}
\def\lefi{\le_{\sf fi}}
\def\cc{\textsf{cc}}
\def\tail{\textsf{tail}}
\def\tt{\textsf{t}}
\def\fim{{\sf fi}}
\def\mal{{\sf mal}}
\def\Comm{{\sf Comm}}
\def\red{\mathop{{\sf red}}}
\def\FIS{{\sf FIS}}

\def\mapright#1{\mathop{\longrightarrow}\limits^{#1}}

\def\istep{\longrightarrow_{{\sf i}}}
\def\restep{\longrightarrow_{{\sf re}}}

\newcommand{\Z}{\mathbb{Z}}

\begin{document}

\title{On finite-index extensions of subgroups of free groups\protect\footnote{%
The first author acknowledges support from Project ASA
(PTDC/MAT/65481/2006) and C.M.U.P., financed by F.C.T. (Portugal)
through the programmes POCTI and POSI, with national and European
Community structural funds.  Both authors acknowledge support from  ESF project
\textsc{AutoMathA}.}%
}

\author{
    Pedro Silva, \small{\url{pvsilva@fc.up.pt}}
    \\
    \small{Centro de Matem\'atica, Universidade do Porto}
    \protect\footnote{%
    Faculdade de Ci\^encias, Universidade do Porto, R. Campo Alegre 687, 
    4169-007 Porto, Portugal}
    \\
    \null
    \\
    Pascal Weil, \small{\url{pascal.weil@labri.fr}}
    \\
    \small{LaBRI, Universit\'e de Bordeaux, CNRS} and \small{CSE 
    Department, IIT Delhi}\protect\footnote{%
    LaBRI, 351 cours de la Lib\'eration, 33400 Talence, France.}
    }

    \date{}

\maketitle

\begin{abstract}
    We study the lattice of finite-index extensions of a given
    finitely generated subgroup $H$ of a free group $F$.  This lattice
    is finite and we give a combinatorial characterization of its
    greatest element, which is the commensurator of $H$.  This
    characterization leads to a fast algorithm to compute the
    commensurator, which is based on a standard algorithm from automata
    theory.  We also give a sub-exponential and super-polynomial upper
    bound for the number of finite-index extensions of $H$, and we
    give a language-theoretic characterization of the lattice of
    finite-index subgroups of $H$.  Finally, we give a polynomial time
    algorithm to compute the malnormal closure of $H$.
\end{abstract}

\bigskip


\medskip


\bigskip



This paper is part of the study of the lattice of finitely generated
subgroups of a free group of finite rank $F$.  Like most of the recent
work on this topic, our paper makes crucial use of the graphical
representation of the subgroups of $F$ introduced in the seminal
papers of Serre (1977 \cite{Serre}) and Stallings (1983
\cite{Stallings}).  This representation not only makes it easier to
form an intuition and to prove properties of subgroups of $F$, but it
also provides a convenient framework to efficiently solve algorithmic
problems and compute invariants concerning these subgroups.

The particular object of study in this paper is the lattice of
extensions of a given finitely generated subgroup $H$ of $F$, and more
specifically the sublattice of finite-index extensions of $H$.  In
this paper, all groups are subgroups of a fixed free group, and the
notion of extension must be understood in this context.

It is elementary to verify that $H$ has only finitely many
finite-index extensions, and it is known that if $K$ and $L$ are
finite-index extensions of $H$, then the subgroup they generate,
namely $\langle K,L\rangle$, has finite index over $H$ as well
(Greenberg's theorem, see \cite{Stallings}).  Therefore $H$ has a
maximum finite-index extension $H_{\fim}$, which is effectively
constructible, and the finite-index extensions of $H$ form a full
convex sublattice of the lattice of subgroups of $F$.

This paper contains a detailed discussion of the lattice of
finite-index extensions of $H$.  Our main contributions are the
following.

We show that the maximum finite-index extension $H_\fim$ of $H$ is the
commensurator of $H$, and we give a combinatorial (graph-theoretic)
characterization of $H_\fim$.  This characterization leads to
efficient algorithms to compute all finite-index extensions of $H$,
and to compute $H_\fim$ -- the latter in time $\O(n\log n)$.  We also
give a rather tight upper bound on the number of finite-index
extensions of $H$: there are at most $\O(\sqrt n\ n^{\frac12\log_2n})$
such extensions, where $n$ is the number of vertices in the graphical
representation of $H$.  Note that this upper bound is sub-exponential
but super-polynomial.

The consideration of the subgroups of the form $H_\fim$, which have no
proper finite-index extensions, leads us to the dual study of the
lattice of finite-index subgroups of a given subgroup, and we give a
combinatorial (language-theoretic) characterization of each such
lattice.

Finally, we use our better understanding of the lattice of extensions
of a subgroup of $F$, to give a polynomial time algorithm to compute
the malnormal closure of a given subgroup.

As we already indicated, we use in a fundamental way the graphical
representation of finitely generated subgroups of $F$, including a
detailed study of the different steps of the computation of this
representation (given a set of generators for the subgroup $H$), whose
study was at the heart of an earlier paper by the authors \cite{SW08}.
It is particularly interesting to see that language-theoretic results
and arguments play an important role in this paper: that is, we
sometimes consider the graphical representation of a subgroup not just
as an edge-labeled graph, but as a finite state automaton.  Such
considerations are present in almost all the results of this paper,
but they become crucial at a rather unexpected juncture: the design of
an efficient algorithm to compute the maximal finite-index extension
$H_\fim$ of $H$.  Indeed, the very low complexity we achieve is due to
the possibility of using a standard automata-theoretic algorithm,
namely the computation of the minimal automaton of a regular language.

Section~\ref{sec: background} summarizes a number of well-known facts
about free groups and the representation of their finitely generated
subgroups, which will be used freely in the sequel (see
\cite{Stallings,Weilsurvey,KM,MVW,SW08} for more details).
Section~\ref{sec: finite index} is the heart of the paper: it starts
with a technical study of the different steps of the algorithm to
compute the graphical representation of a given subgroup, and a
description of those steps which preserve finite-index
(Section~\ref{preserve fi}).  These technical results are then used to
characterize the maximal finite-index extension $H_\fim$
(Section~\ref{lattice fi}), to relate the computation of $H_\fim$ and
the minimization of certain finite-state automata (Section~\ref{sec:
computing}), to evaluate the maximal number of finite-index extensions
of a given subgroup (Section~\ref{sec: counting}), and to describe an
invariant of the lattice of finite-index subgroups of a given subgroup
(Section~\ref{sec: FIS}).

Finally, we apply the same ideas in Section~\ref{sec: malnormal}, to
study the malnormal closure of a subgroup, and to show that it can be
computed in polynomial time.

\section{Subgroups of free groups and Stallings graphs}\label{sec: 
background}

Let $F$ be a finitely generated free group and let $A = \{a_1, \ldots,
a_r\}$ be a fixed basis of $F$.  Let $\bar A = \{\bar a_1, \ldots,
\bar a_r\}$ be a disjoint copy of $A$ and let $\tilde A = A \cup \bar
A$: as usual, we extend the map $a\mapsto \bar a$ from the set $A$ to
all words by letting $\bar{\bar a} = a$ if $a\in A$ and $\overline{ua}
= \bar a\bar u$ if $a\in \tilde A$ and $u\in \tilde A^*$.  As usual
again, the elements of $F$ are identified with the \textit{reduced
words} over the alphabet $\tilde A$, that is, the words that do not
contain a sequence of the form $a\bar a$ ($a\in \tilde A$).  If $u\in
\tilde A^*$ is an arbitrary word, we denote by $\red(u)$ the
corresponding reduced word, that is, the word obtained from $u$ by
repeatedly deleting all sequences of the form $a\bar a$ ($a\in \tilde
A$).

A reduced word $u \in F$ is \textit{cyclically reduced} if $u$ cannot
be written as $u = av\bar a$ with $a\in \tilde A$ and $v \in F$.
Every reduced word $u$ can be factored in a unique way in the form $u
= xy\bar x$, with $y$ cyclically reduced.

If $H$ is a subgroup of $F$, an \textit{extension} of $H$ is any
subgroup $G$ containing $H$ and we write $H \le G$.  If $H$ is
finitely generated, we also write $H \lefg G$.  If $H$ has finite
index in $G$, we say that $G$ is a \textit{finite-index extension} of
$H$ and we write $H \lefi G$. Finally, we write $H \leff G$ if $H$ is 
a free factor of $G$.

\subsection{The graphical representation of a subgroup}

It is well known (since Serre's and Stalling's fundamental work
\cite{Serre,Stallings}) that every finitely generated subgroup $H
\lefg F$ admits a unique graphical representation of the form
$\calA(H) = (\Gamma(H),1)$, where $\Gamma(H)$ is a finite directed
graph with $A$-labeled edges and 1 is a designated vertex of
$\Gamma(H)$, subject to the combinatorial conditions below.  Here, a
graph is a pair $(V,E)$ where $V$ is the set of \textit{vertices} and
$E \subseteq V \times A \times V$ is the set of \textit{edges}; the
\textit{in-degree} (resp.  \textit{out-degree}) of a vertex $v\in V$
is the number of edges in $E$ of the form $(v',a,v)$ (resp.
$(v,a,v')$); and the \textit{degree} of $v$ is the sum of its in- and
out-degree.  Every pair $\calA(H)$ satisfies the following:

- the (underlying undirected) graph is connected;

- for each $a\in A$, every vertex is the source (resp.  the target) of
at most one $a$-labeled edge;

- and every vertex, except possibly 1, has degree at least 2.

\smallskip\noindent Moreover, every pair $(\Gamma,1)$ with these
properties is said to be \textit{admissible}, and it is the
representation of a finitely generated subgroup of $F$.  In addition,
given a finite set of generators of $H$, the representation of $H$ is
effectively computable.  We refer the reader to
\cite{Stallings,KM,MVW,Weilsurvey,SW08} for some of the literature on
this construction and its many applications, and to Section~\ref{sec:
fi extensions} below on the construction of $\calA(H)$.

We sometimes like to view the $A$-labeled graph $\Gamma(H)$ as a
transition system over the alphabet $\tilde A$: if $p,q$ are vertices of
$\Gamma(H)$, $a\in A$ and $(p,a,q)$ is an edge of $\Gamma(H)$, we
say that $a$ labels a path from $p$ to $q$ and $\bar a$ labels a path
from $q$ to $p$, written $p \mapright a q$ and $q \mapright{\bar a}
p$.  If $u \in \tilde A^*$ is a word (reduced or not) and $u = va$
($a\in \tilde A$), we say that $p \mapright u q$ (\textit{$u$ labels a
path from $p$ to $q$}) if $p \mapright v p' \mapright a q$ for some
vertex $p'$.  In particular, a reduced word is in $H$ if and only if
it labels a loop at vertex 1.  Moreover, if we have a path $p
\mapright u q$, then we also have a path $p \mapright{\red(u)} q$.

If $H \lefg G \lefg F$, then there is a \textit{homomorphism} from
$\calA(H)$ into $\calA(G)$, that is, a map $\phi$ from the vertex set 
of $\Gamma(H)$ to the vertex set 
of $\Gamma(G)$ such that

- $\phi(1) = 1$ and

- if $p \mapright a q$ in $\Gamma(H)$ ($p,q$ vertices, $a\in \tilde 
A$), then $\phi(p) \mapright a \phi(q)$ in $\Gamma(G)$.

\smallskip\noindent It is not difficult to verify that this morphism, if it
exists, is unique, and we denote it by $\phi_{H}^{G}$.
It is well known (see \cite{Serre,Stallings,KM,MVW}) that if
$\phi_H^G$ is one-to-one, then $H \leff G$.

Finally, we say that the homomorphism $\phi_H^G$ is a \textit{cover},
if it satisfies

- if $p,q$ are vertices of $\Gamma(H)$, $a\in \tilde A$ and
$\phi_H^G(p) \mapright a \phi_H^G(q)$ in $\Gamma(G)$, then $p
\mapright a q'$ in $\Gamma(H)$ for some vertex $q'$ such that
$\phi_H^G(q') = \phi_H^G(q)$.  In that case, all sets of the form
$(\phi_H^G)\inv(q)$ ($q$ a vertex of $\Gamma(G)$) have the same
cardinality.

Covers have the following property, which we will use freely in the 
sequel.

\begin{lemma}
    If $\phi\colon \calA(H) \to \calA(G)$ is a cover, $p$ is a vertex
    of $\calA(H)$ and $u \in F$ labels a loop of $\calA(G)$ at
    $\phi(p)$, then $\red(u^m)$ labels a loop of $\calA(H)$ at $p$ for
    some integer $m > 1$.
\end{lemma}

\subsection{Covers, cyclically reduced subgroups and finite-index extensions}

Let us say that $H$ is \textit{cyclically reduced} (with respect to
the basis $A$) if every vertex of $\Gamma(H)$ has degree at least
equal to 2.  If $H$ is not cyclically reduced, then the designated
vertex 1 of $\Gamma(H)$ has degree 1 and $\Gamma(H)$ consists of two
parts: $\tail(\Gamma(H))$, which contains the designated vertex 1 and
all degree 2 vertices that can be connected to vertex 1 through other
degree 2 vertices; and the rest of $\Gamma(H)$, which is called the
\textit{core} of $\Gamma(H)$, written $\cc(\Gamma(H))$.  We let
$\tt_H(1)$ be the shortest word which labels a path from 1 to a vertex
in $\cc(\Gamma(H))$ and let $\tau_H(1)$ be the vertex of
$\cc(\Gamma(H))$ thus reached (if $H$ is cyclically reduced, then
$\tt_H(1)$ is the empty word and $\tau_H(1) = 1$).  We write $\tt(1)$
and $\tau(1)$ if the subgroup $H$ is clear from the context.

The tail and the core of $\Gamma(H)$ have intrinsic characterizations.
The characterization of the core is well-known (see \cite[Exercise
7.3(a)]{Stallings}) and that of the tail is an elementary consequence.

\begin{remark}\label{loop and tail}
    Let $H \lefg F$ and let $p$ be a vertex of $\Gamma(H)$. Then 
    $p$ is a vertex of $\cc(\Gamma(H))$ if and only if some 
    cyclically reduced word $u$ labels a path from $p$ to $p$.
\end{remark}

\proof
By definition, if $H$ is not cyclically reduced, then
$\tail(\Gamma(H))$ consists of a single path from vertex 1 to vertex
$\tt(1)$ (excluding the latter vertex): it is therefore elementary to 
verify that no non-empty cyclically reduced word labels a loop at a 
vertex in $\tail(\Gamma(H))$.

Let now $p$ be a vertex in $\cc(\Gamma(H))$: then $p$ has degree at
least 2, and if it has degree exactly 2, then neither of the two edges
adjacent to it leads to a vertex in $\tail(\Gamma(H))$.  Therefore,
one can find distinct letters $a,a'\in \tilde A$ such that $p
\mapright a q$ and $p \mapright{a'} q'$, with $q$ and $q'$ in
$\cc(\Gamma(H))$ as well.  Iterating this reasoning, one can show that
there exist arbitrarily long paths within $\cc(\Gamma(H))$, starting
from $p$ and labeled by reduced words of the form $au$ and $a'u'$.
Since $\Gamma(H)$ is finite, vertices are repeated along these paths,
and we consider the earliest such repetition after the initial $p$.
If $p$ itself is the first repeated vertex along the path labeled
$au$, we have a loop $p \mapright{au} p$ such that $au$ is cyclically
reduced, and we are done.  The situation is similar if $p$ is the
first repeated vertex along the path labeled $a'u'$.  Otherwise, let
$r$ and $r'$ be the first repeated vertices along the two paths.  Then
$r \ne p$, $r' \ne p$, and $\cc(\Gamma(H))$ has paths of the form $p
\mapright{au} r$, $p \mapright{a'u'} r'$, $r \mapright v r$ and $r'
\mapright{v'} r'$ such that $auv\bar u\bar a, a'u'v'\bar u'\bar a'$
are reduced.  Then the word $auv\bar u\bar aa'u'v'\bar u'\bar a'$ is
cyclically reduced, and it labels a loop at $p$ in $\Gamma(H)$.  This
concludes the proof.
\eop

\begin{remark}\label{tail only}
    Let $H \lefg F$.  Then $\tt(1)$ is the maximum common prefix of
    the non-trivial elements of $H$.
\end{remark}

\proof
Since every non-trivial element of $H$ is the label of a loop at 1 in
$\Gamma(H)$, it is clear that $\tt(1)$ is a common prefix to all these
words.

By Remark~\ref{loop and tail}, there exists a cyclically reduced word
$u$ labeling a loop at $\tau(1)$.  Then both
$\tt(1)u\overline{\tt(1)}$ and $\tt(1)\bar u\overline{\tt(1)}$ are
reduced words in $H$, and their maximum common prefix is $\tt(1)$.
This concludes the proof.
\eop

We can now state the following extension of the classical 
characterization of finite-index extensions of cyclically reduced 
subgroups in terms of covers.

\begin{proposition}\label{prop: fi cover}
    Let $H\lefg G \lefg F$.  Then $\tt_G(1)$ is a prefix of $\tt(H)$.
    Moreover, $H \lefi G$ if and only if $\tt_H(1) = \tt_G(1)$ and the
    restriction of $\phi^G_H$ is a cover from $(\cc(\Gamma(H)),
    \tau_H(1))$ onto $(\cc(\Gamma(G)), \tau_G(1))$.  If that is the
    case, the index of $H$ in $G$ is the common cardinality of the
    subsets ${\phi^G_H}\inv(q)$ ($q$ a vertex of $\cc(\Gamma(G))$).
\end{proposition}

\proof
Let $\phi = \phi^G_H$.
If $u$ is cyclically reduced and labels a loop at a vertex $p$ of
$\cc(\Gamma(H))$, then $u$ labels a loop at $\phi(p)$ in $\Gamma(G)$,
and that vertex is in the core of $\Gamma(G)$ by Fact~\ref{loop and
tail}.  It follows that $\tt_H(1)$ labels a path from the origin in
$\Gamma(G)$ to a vertex in $\cc(\Gamma(G))$.  In particular,
$\tt_G(1)$ is a prefix of $\tt_H(1)$.

If $\tt_G(1)$ is a proper prefix of $\tt_H(1)$, we have $\tt_H(1) =
\tt_G(1)at$ for some $a\in \tilde A$ and $t\in F$.  Since $\tau_G(1)$
is in $\cc(\Gamma(G))$, there exist a cyclically reduced word of the
form $bu$, with first letter $b\ne a$, which labels a loop at
$\tau_G(1)$ in $\cc(\Gamma(G))$.  Then the words
$\tt_G(1)(bu)^n\overline{\tt_G(1)}$ are all reduced, and the cosets
$H\tt_G(1)(bu)^n\overline{\tt_G(1)}$ are all in $G$.  Moreover, these
cosets are pairwise disjoint since $H$ contains no reduced word of the
form $\tt_G(1)(bu)^d\overline{\tt_G(1)}$, $d\ne 0$.  Thus, if $H\lefi
G$, then $\tt_H(1) = \tt_G(1)$.  It follows immediately that $\phi$
maps core vertices to core vertices and tail vertices to tail
vertices.

Let us now assume that $\Gamma(H)$ and $\Gamma(G)$ have the same
tails, and let us denote by $\tt(1)$ the word $\tt_H(1) = \tt_G(1)$.
Then $H' = \overline{\tt(1)}H\tt(1)$ and $G' =
\overline{\tt(1)}G\tt(1)$ are cyclically reduced, and $H\lefi G$ if
and only if $H' \lefi G'$.  Thus we may now assume that $G$ and $H$
are cyclically reduced.  If $\phi$ is not a cover, there exists a
vertex $p$ of $\Gamma(H)$ such that $\Gamma(G)$ has a loop at
$\phi(p)$ labeled by a cyclically reduced word $bu$ ($b\in \tilde A$)
and $\Gamma(H)$ has no $b$-labeled edge out of $p$.  Let $v$ label a
path from 1 to $p$ in $\Gamma(H)$ (and hence in $\Gamma(G)$).  By the
same reasoning as above, the cosets $Hv(bu)^n\bar v$ are pairwise
distinct, and contained in $G$. Thus, if $H \lefi G$, then $\phi$ is 
a cover.

The converse is verified as follows: if $\phi$ is a cover, let 
$u_1,\ldots,u_d$ be reduced words labeling paths in $\Gamma(H)$ from 
$1$ to the elements $1 = p_1,\ldots,p_d$ of $\phi\inv(1)$. If $g\in 
G$, then $g$ labels a loop at 1 in $\Gamma(G)$, and since $\phi$ is a 
cover, $g$ labels a path in $\Gamma(H)$ from 1 to $p_i$ for some $i$. 
Therefore $g\in Hu_i$: thus $G$ is the union of finitely many 
$H$-cosets.
\eop

\begin{corollary}
    The extensions, and the finite-index subgroups of a cyclically
    reduced subgroup are cyclically reduced as well.
\end{corollary}

\section{Finite-index extensions of a subgroup $H$}\label{sec: finite
index}

It follows from the characterization of finite-index extensions by
covers, that if $H \lefi G$, then $\phi_H^G$
is onto\footnote{%
The converse is not true, see \cite{MVW} for a detailed study of the 
extensions $H \lefg G$ such that $\phi_H^G$ is onto.%
}.  Therefore $H$ has only a finite number of finite-index extensions,
and that number can be bounded above by the number of binary relations
on $\Gamma(H)$: if that graph has $n$ vertices, then $H$ has at most
$2^{n^2}$ finite-index extensions.  We give a better upper bound in
Section~\ref{sec: counting}.

Moreover, the collection of finite-index extensions of $H$ is
effectively computable.  In addition, it is elementary to use these
graphical representations to show that the join of two finite-index
extensions of $H$ is again a finite-index extension (see Stalling's
proof of Greenberg's theorem \cite{Stallings}).  It follows that if $H
\lefg F$, then $H$ admits an effectively computable maximum
finite-index extension $H_\fim$.  The results of Section~\ref{sec:
computing} below yield an efficient algorithm to compute the set of
finite-index extensions of a given subgroup $H$, and its maximum
finite-index extension $H_\fim$.

\begin{remark}\label{remark index}
    We just observed that every finitely generated subgroup of $F$ has
    a finite number of finite-index extensions.  However, it usually
    has infinitely many finite-index subgroups.  More precisely, every
    non-trivial subgroup $H \lefg F$ admits a finite-index subgroup of
    index $r$ for each $r \ge 1$.  Indeed, let $\calA(H) =
    (\Gamma(H),1)$ and let $Q$ be the vertex set of $\Gamma(H)$.
    Define $\Gamma_r$ be the $A$-labeled graph with vertex set $Q
    \times \{1,\ldots,r\}$ and with the following edge set: for each
    edge $(p,a,q)$ of $\Gamma(H)$, there is an edge
    $((p,i),a,(q,i+1))$ for each $1\le i < r$ and an edge
    $((p,r),a,(q,1))$.  Then $\Gamma_r$ is an admissible graph, the
    map $\pi\colon (p,i) \mapsto p$ defines a cover from $\Gamma_r$ to
    $\Gamma(H)$, and if $H_r$ is the subgroup represented by
    $(\Gamma_r, (1,1))$, then $H_r$ has index $r$ in $H$.
\end{remark}

\subsection{i-steps and finite-index extensions}\label{sec: fi extensions}

Let $H \le G$ be finitely generated subgroups of $F$ and let $g_1,
\ldots, g_n$ be reduced words such that $G = \langle H, g_1, \ldots,
g_n\rangle$.  Let $G_0 = H$ and let $G_i = \langle G_{i-1},
g_i\rangle$ ($1 \le i \le n$).  We may of course assume that $g_i
\not\in G_{i-1}$, so $G_{i-1} \ne G_i$.

Then $\calA(G_i)$ is obtained from $\calA(G_{i-1})$ by, first, adding
sufficiently many new vertices and edges to create a new path from
vertex 1 to itself, labeled by $g_i$; and second, by 
\textit{reducing}\footnote{%
This reduction operation is the iteration of Stallings's
\textit{folding} operation \cite{Stallings}; our terminology
emphasizes the fact that this is a generalization of the reduction of
a word (the iterated process of deleting factors of the form $a\bar
a$, $a\in \tilde A$).}
the resulting graph, that is, repeatedly identifying vertices $p$ and
$p'$ such that $q \mapright a p$ and $q \mapright a p'$ for some
vertex $q$ and some letter $a\in \tilde A$, see for instance
\cite{Stallings,Weilsurvey,KM,SW08}.  Depending on the length of
prefixes of $g_i$ and $\bar g_i$ that can be read from vertex $1$ in
$\Gamma(G_{i-1})$, this procedure amounts to one of the two
following moves:

- a \textit{reduced expansion}, or \textit{re-step} (we write
$\calA(G_{i-1}) \restep^{(p,w,q)} \calA(G_i)$), that is, we add a new
path labeled by a factor $w$ of $g_i$, from some vertex $p$ to some
vertex $q$ of $\Gamma(G_{i-1})$ in such a way that the resulting
graph is admissible (needs no reduction);

- or an \textit{i-step} (we write $\calA(G_{i-1}) \istep^{p=q}
\calA(G_i)$), that is, we identify a pair of vertices $(p,q)$ of
$\Gamma(G_{i-1})$, and we reduce the resulting graph.

\begin{remark}
    Let us comment on these steps, with reference to Stallings's
    algorithm \cite{Stallings}.  If $H = \langle
    g_1,\ldots,g_n\rangle$, Stallings produces $\calA(H)$ by reducing
    (folding) a bouquet of $n$ circles, labeled $g_1,\ldots, g_n$
    respectively.  For our purpose, we decompose this operation in $n$
    steps, adding one generator at a time and producing successively
    the $\calA(\langle h_1,\ldots,h_i\rangle)$ ($1\le i\le n$).  Each
    of these steps is either an re-step or an i-step.
\end{remark}    

We refer the readers to \cite[Section 2]{SW08} for a detailed analysis
of these moves and we record the following observation.

\begin{lemma}\label{first lemma i-steps}
    Let $G,H \lefg F$. If $H \lefi G$, then only i-steps are 
    involved in the transformation from $\Gamma(H)$ to $\Gamma(G)$.
\end{lemma}

\proof
Let $H = G_0 \le G_1 \le \ldots \le G_n = G$ be as in the above
discussion.  Note that $H \lefi G$ if and only if $G_{i-1} \lefi G_i$
for each $1 \le i \le n$.  If $\calA(G_{i-1}) \restep \calA(G_i)$,
then the homomorphism $\phi_{G_{i-1}}^{G_i}$ is one-to-one, so
$G_{i-1} \leff G_i$ and in particular, $G_i$ is not a finite-index
extension of $G_{i-1}$.
\eop

\subsection{Which i-steps yield finite-index 
extensions?}\label{preserve fi}

If $p$ is a vertex of $\cc(\Gamma(H))$, we let $\widetilde L_p(H)$ be
the language accepted by $\cc(\Gamma(H))$, seen as a finite state
automaton with initial state $p$ and all states final: that is, the
set of (possibly non-reduced) words in $\tilde A^*$ that label a path
in $\cc(\Gamma(H))$ starting at vertex $p$.  Let then $L_p(H)$ be the
set of reduced words in $\widetilde L_p(H)$ --- which is also the set
of all $\red(u)$ ($u \in \widetilde L_p(H)$), and also the set of
prefixes of words in the subgroup represented by the pair
$(\cc(\Gamma(H)),p)$.  Let us first record the following elementary
remark.

\begin{remark}\label{fact: Lp Lq}
    If $p \mapright u q$ is a path in $\cc(\Gamma(H))$, then $L_q(H) =
    \{\red(\bar ux) \mid x \in L_p(H)\}$.
\end{remark}    

We now refine the result of Lemma~\ref{first lemma i-steps}.

\begin{proposition}\label{cns finite index}
    Let $H\lefg F$ and let $p,q$ be distinct vertices in $\Gamma(H)$.
    Let $G$ be the subgroup of $F$ such that $\calA(H) \istep^{p=q}
    \calA(G)$.  Then $H \lefi G$ if and only if $p,q$ are in
    $\cc(\Gamma(H))$ and $L_p(H) = L_q(H)$, if and only if $p,q$ are
    in $\cc(\Gamma(H))$ and $\widetilde L_p(H) = \widetilde L_q(H)$.
\end{proposition}

\proof
Let us first assume that $H \lefi G$ and let $\phi = \phi^G_H$.  By
Proposition~\ref{prop: fi cover}, $\phi$ is a bijection from
$\tail(H)$ onto $\tail(G)$ and, since $\phi(p) = \phi(q)$, the
vertices $p$ and $q$ must both be in $\cc(\Gamma(H))$.

\smallskip

\textbf{If $p,q \in \cc(\Gamma(H))$ and $\widetilde L_p(H) \ne
\widetilde L_q(H)$}, we consider (without loss of generality) a word
$u \in \widetilde L_p(H) \setminus \widetilde L_q(H)$, with minimum
length, say $u = va$ with $a\in \tilde A$.  By definition, there exist
paths $p \mapright v p' \mapright a p''$ and $q \mapright v q'$, but
no path $q' \mapright a q''$ in $\cc(\Gamma(H))$.  Observe that
$\calA(H) \istep^{p' = q'} \calA(G)$.  If there is a path $q'
\mapright a q''$ in $\Gamma(H)$, then we also have $\calA(H)
\istep^{p'' = q''} \calA(G)$ and since $p'' \in \cc(\Gamma(H))$ and
$q'' \in \tail(\Gamma(H))$, we conclude to a contradiction by
Proposition~\ref{prop: fi cover}.  We now assume that there is no path
$q' \mapright a q''$ in $\Gamma(H)$.

We claim that there exists $w\in F$ such that $aw$ is cyclically
reduced and $p' \mapright{aw} q'$ in $\cc(\Gamma(H))$.  Let indeed
$p'' \mapright z q'$ be a path of minimal length in $\cc(\Gamma(H))$
(there exists one by connectedness).  Since $p'' \in \cc(\Gamma(H))$,
there exists a path $p'' \mapright b r$ for some $b \in \tilde A$,
$b\ne \bar a$, and as in the proof of Remark~\ref{loop and tail}, there
exists a reduced word of the form $bt$ labeling a loop at $p''$.  Let
$w = \red(btz)$: then we have a path $p'' \mapright w q'$.  By
minimality of the length of $z$, $\bar t\,\bar b$ is not a prefix of
$z$, so $w$ starts with letter $b$, and hence $aw$ is reduced.  In
fact, $aw$ is cyclically reduced since there is no path $q' \mapright
a q''$.

Let $1 \mapright t q'$ be a path in $\Gamma(H)$.  Then $\red(taw\bar
t) \in G$, and hence there exists $m > 1$ such that
$\red(t(aw)^m\bar t) \in H$.  Again, since there is no path $q'
\mapright a q''$, the word $ta$ is reduced. By replacing $m$ by a 
sufficiently large multiple, we find that $ta$ is a prefix of 
$\red(t(aw)^m\bar t)$, and hence that $ta$ labels a path from 1 in 
$\Gamma(H)$: this contradicts the absence of a path $q' \mapright a q''$.

\smallskip

Thus we have proved that, if $H \lefi G$, then $p,q \in
\cc(\Gamma(H))$ and $\widetilde L_p(H) = \widetilde L_q(H)$.  The
latter condition immediately implies that $L_p(H) = L_q(H)$.
\textbf{We now assume that $p,q \in \cc(\Gamma(H))$ and $L_p(H) =
L_q(H)$, and we show that $H \lefi G$}.  We first establish a
technical fact.

\begin{lemma}\label{special fact}
    Let $r_i \mapright{z_i} s_{i+1}$ ($z_i \in F$, $0 \le i \le k$) be
    paths in $\cc(\Gamma(H))$, such that $r_i, s_i \in \{p,q\}$ for
    each $1\le i \le k$.  Then there exists a path $r_0
    \mapright{\red(z_0\cdots z_k)} t$ in $\cc(\Gamma(H))$.
\end{lemma}

\proof
The proof is by induction on $k$, and is trivial for $k = 0$.  If $k >
0$, then there is a path $r_1 \mapright{\red(z_1\cdots z_k)} t$ in
$\cc(\Gamma(H))$.  Since $L_p(H) = L_q(H)$, there is also a path $s_1
\mapright{\red(z_1\cdots z_k)} t'$ for some $t' \in \cc(\Gamma(H))$,
and therefore a path $r_0 \mapright{\red(z_0\cdots z_k)} t'$ as
required.
\eop

We want to show that $G$ has finitely many $H$-cosets.  Let $u \in G$:
then $u$ labels a loop at 1 in $\Gamma(G)$. Let $\calB$ be the 
automaton obtained from $\Gamma(H)$ by identifying vertices $p$ and 
$q$, but without performing any reduction. Then $\Gamma(G)$ is the 
result of the reduction of $\calB$. In particular (say, in view of 
\cite[Fact 1.4]{SW08}), $u = \red(v)$ for some word $v \in \tilde 
A^*$ labeling a loop at 1 in $\calB$. By definition of $\calB$, the 
word $v$ factors as
$v = v_0\cdots v_k$, in such a way that $\Gamma(H)$ has paths of the 
form $1 \mapright{v_0} s_1$, $r_i \mapright{v_i} s_{i+1}$ ($1
\le i < k$) and $r_k \mapright{v_k} 1$, and the vertices $r_1, s_1, 
\cdots, r_k, s_k$ are all equal to $p$ or $q$. As observed in 
Section~\ref{sec: 
background},  $\Gamma(H)$ also has paths
$$1 \mapright{\red(v_0)} s_1,\ r_i \mapright{\red(v_i)} s_{i+1}\ (1
\le i < k)\textrm{ and }r_k \mapright{\red(v_k)} 1.$$
In particular, we have $\red(v_0) = \tt(1)w_0$ and $\red(v_k) =
w_k\bar\tt(1)$ for some $w_0, w_k \in F$, and there are paths $1
\mapright{\tt(1)} \tau(1) \mapright{w_0} s_1$ and $r_k \mapright{w_k}
\tau(1) \mapright{\bar\tt(1)} 1$.  Note that the paths $\tau(1)
\mapright{w_0} s_1$, $r_i \mapright{\red(v_i)} s_{i+1}$ ($1 \le i <
k$) and $r_k \mapright{w_k} \tau(1)$ are set entirely within
$\cc(\Gamma(H))$, since no reduced word-labeled path between vertices
in $\cc(\Gamma(H))$ can visit a vertex in $\tail(\Gamma(H))$.

By Lemma~\ref{special fact}, there exists a path $\tau(1)
\mapright{\red(w_0v_1\cdots v_{k-1}w_k)} t$ for some vertex $t$ in
$\cc(\Gamma(H))$.  Let $h$ be a shortest-length word such that $t
\mapright h \tau(1)$ in $\cc(\Gamma(H))$.  Then $z =
\tt(1)\red(w_0v_1\cdots v_{k-1}w_k)\ h\ \bar\tt(1)$ labels a loop at
vertex 1 in $\Gamma(H)$, so $\red(z) \in H$.  By construction, we have
$u = \red(\tt(1)w_0v_1\cdots v_{k-1}w_k\bar\tt(1))$, so $u \in
H\red(\tt(1)\ \bar h\ \bar\tt(1))$.  Since $h$ was chosen to be a geodesic
in $\cc(\Gamma(H))$, it can take only finitely many values, and this
completes the proof that $G$ has finitely many $H$-cosets.
\eop

We note the following consequence of this proof.

\begin{corollary}\label{fact tilde or not}
    Let $H \lefg F$ and let $p,q \in \cc(\Gamma(H))$. Then $L_p(H) = 
    L_q(H)$ if and only if $\widetilde L_p(H) = \widetilde L_q(H)$.
\end{corollary}    

\subsection{The lattice of finite-index extensions of $H$}\label{lattice fi}

We further refine Proposition~\ref{cns finite index} as follows: we
consider an extension $H\lefi G$ and a pair $(r,s)$ of vertices of
$\Gamma(H)$, whose identification yields a finite-index extension of
$H$. Then we show that identifying the
vertices of $\Gamma(G)$ corresponding to $r$ and $s$, also yields a
finite-index extension of $G$.

\begin{lemma}\label{preserve Lr}
    Let $H \lefi G \lefg F$ and let $\phi = \phi^G_H$.  Let
    $p$ be a vertex of $\cc(\Gamma(H))$.  Then $\widetilde L_p(H) =
    \widetilde L_{\phi(p)}(G)$ and $L_p(H) = L_{\phi(p)}(G)$.
\end{lemma}

\proof
If $p \mapright u r$ in $\cc(\Gamma(H))$, then the $\phi$-image of
this path is a path $\phi(p) \mapright u \phi(r)$, which is entirely
contained in $\cc(\Gamma(G))$ by Proposition~\ref{prop: fi cover}.  In
particular, $\widetilde L_p(H) \subseteq \widetilde L_{\phi(p)}(G)$.

Conversely, suppose that $\phi(p) \mapright u r'$ is a path in
$\cc(\Gamma(G))$.  Since $\phi$ is a cover from $\cc(\Gamma(H))$ onto
$\cc(\Gamma(H))$, $u$ labels some path $p \mapright u r$ in
$\cc(\Gamma(H))$, and hence $u \in \widetilde L_p(H)$.  Thus
$\widetilde L_p(H) = \widetilde L_{\phi(p)}(G)$.
\eop

Together with Proposition~\ref{cns finite index}, Lemma~\ref{preserve 
Lr} immediately implies the following statements.

\begin{corollary}\label{charact fi extensions}
    Let $H \lefg F$.
    
    \begin{enumerate}
	\item Let $H \lefi G$ and let $\phi = \phi^G_H$.  If $p,q$ are
	vertices of $\Gamma(H)$, $\calA(H) \istep^{p = q} \calA(K)$
	and $\calA(G) \istep^{\phi(p) = \phi(q)} \calA(K')$, then $H
	\lefi K$ if and only if $G \lefi K'$.
	
	\item $H \lefi G$ if and only if $\calA(G)$ is obtained from 
	$\calA(H)$ by identifying some pairs of vertices $(p,q)$ in 
	$\cc(\Gamma(H))$ such that $L_p(H) = L_q(H)$, and then 
	reducing the resulting graph.
    \end{enumerate}	
\end{corollary}

The identification of all pairs of vertices $(p,q)$ such that $L_p(H)
= L_q(H)$ yields the minimum quotient of $\calA(H)$ and hence the
maximum finite-index extension $H_\fim$ of $H$.  In addition, we find
that $H_\fim$ is exactly the \textit{commensurator} of $H$ (the set
$\Comm_F(H)$ of all elements $g\in F$ such that $H \cap H^g$ has finite
index in both $H$ and $H^g$), a fact that can also be deduced from
\cite[Lemma 8.7]{KM}.

\begin{theorem}\label{maximal fi extension}
    Let $H\lefg F$ and let $H_\fim$ be its maximum finite-index
    extension.
    
    \begin{enumerate}
	\item $\calA(H_\fim)$ is obtained from $\calA(H)$ by
	identifying all pairs of vertices $p,q$ of $\cc(\Gamma(H))$
	such that $L_p(H) = L_q(H)$.  No reduction is necessary.
	
	\item $H_\fim = \Comm_F(H)$.
    \end{enumerate}	
\end{theorem}

\proof
In view of Corollary~\ref{charact fi extensions}, $\calA(H_\fim)$ is
obtained from $\calA(H)$ by identifying all pairs of vertices $p,q$ of
$\cc(\Gamma(H))$ such that $L_p(H) = L_q(H)$, and then by reducing the
resulting graph $\calB$.  If $p \mapright a r$ and $q \mapright a s$
($a\in \tilde A$) are paths in $\cc(\Gamma(H))$ and if $L_p(H) =
L_q(H)$, then $L_r(H) = L_s(H)$ by Remark~\ref{fact: Lp Lq}.  Thus
$\calB$ is already reduced, which concludes the proof of the first
statement.

The fact that $\Comm(H)$ is a subgroup and a finite-index extension of
$H$ is proved, for instance, in \cite[Prop.  8.9]{KM}.  Conversely,
suppose that $H \lefi G$ and $g \in G$.  Since conjugation by $g$ is
an automorphism of $G$, we have $H^g \lefi G$.  Now the intersection
of finite-index subgroups, again has finite index, so $H \cap H^g
\lefi G$ and hence $H \cap H^g \lefi H$ and $H \cap H^g \lefi H^g$.
Thus $g \in \Comm_F(H)$, which concludes the proof.
\eop

\subsection{Computing finite-index extensions}\label{sec: computing}

Recall the notion of minimization of a deterministic finite-state
automaton (see \cite{Kozen} for instance).  Let $\calB = (Q,i,E,T)$ be
such an automaton, over the alphabet $B$, with $Q$ the finite set of
states, $i \in Q$ the initial state, $E \subseteq Q \times B \times Q$
the set of transitions and $T \subseteq Q$ the set of accepting
states, and let $L$ be the language accepted by $\calB$, that is, the
set of words in $B^*$ that label a path from $p$ to a state in $T$.
Then the minimal automaton of $L$ is obtained by identifying the pairs
of states $(p,q)$ such that the automata $(Q,p,E,T)$ and $(Q,q,E,T)$
accept the same language.

In our situation, the alphabet is $\tilde A$ and $\widetilde L_p(H)$
is the language accepted by the automaton $\calB_p$, whose states and
transitions are given by $\cc(\Gamma(H))$, with initial state $p$ and
all states final.  Therefore Corollary~\ref{charact fi extensions} and
Theorem~\ref{maximal fi extension} show that the identification of two
vertices $p,q \in \cc(\Gamma(H))$ yields a finite-index extension if
and only if $p$ and $q$ are identified when minimizing
$\calB_{\tau(1)}$.  Moreover, $\cc(\Gamma(H_\fim))$ is given by the
states and transitions of the minimal automaton of $\widetilde
L_{\tau(1)}(H)$.

The classical Hopcroft algorithm (see \cite{Kozen}) minimizes an
$n$-state automaton in time $\O(n\log n)$, so we have the following
result.

\begin{proposition}\label{prop compute Hfi}
    Let $H \lefg F$, and let $n$ be the number of vertices of
    $\Gamma(H)$.
    
    \begin{itemize}
	\item $\cc(\Gamma(H_\fim))$ is obtained by minimizing the 
	automaton given by the vertices and edges of 
	$\cc(\Gamma(H))$, with all states final (the initial state 
	does not matter in that situation).
	
	\item One can compute $\Gamma(H_\fim)$ in time $\O(n\log n)$.
	
	\item One can decide in time $\O(n\log n)$ whether identifying
	a given set of pairs of vertices of $\Gamma(H)$ will produce a
	finite-index extension of $H$.
    \end{itemize}	
\end{proposition}

\begin{remark}
    It may be that for the particular automata at hand (over a
    symmetrized alphabet, with all states final), the complexity of
    Hopcroft's algorithm might be better than $\O(n\log n)$, even
    linear.  It has also been observed that in many instances,
    Myhill's automata minimization algorithm exhibits a better
    performance than Hopcroft's, in spite of a $\O(n^2)$ worst-case
    complexity.  Brzozowski's algorithm \cite{Brz} also performs
    remarkably well in practice \cite{ChZi}.
\end{remark}    

\subsection{Counting finite-index extensions}\label{sec: counting}

Recall that, if $\Gamma$ is an $A$-labeled graph, the \textit{product}
$\Gamma \times_A \Gamma$ (also called the \textit{fiber product}, or
the \textit{pull-back}, of two copies of $\Gamma$) is the $A$-labeled
graph whose vertex set is the set of pairs $(p,q)$ of vertices of
$\Gamma$ and whose edges are the triples $((p,q),a,(p',q'))$ such that
$(p,a,q)$ and $(p',a,q')$ are edges of $\Gamma$.  This graph is not
admissible, nor even connected in general (the vertices of the form
$(p,p)$ form a connected component that is isomorphic to $\Gamma$).
Note that there is a $u$-labeled path in $\Gamma \times_A \Gamma$ from
$(p,q)$ to $(p',q')$, if and only if $\Gamma$ has paths $p \mapright u
p'$ and $q \mapright u q'$.

If $p,q$ are vertices of $\cc(\Gamma(H))$, we let $p \sim q$ if and
only if $L_p(H) = L_q(H)$.

\begin{proposition}\label{LpLq covers}
    Let $H \lefg F$.
    
    \begin{itemize}
	\item The relation $\sim$ is a union of connected components
	of $\cc(\Gamma(H)) \times_A \cc(\Gamma(H))$.
	
	\item Let $p,q$ be vertices of $\cc(\Gamma(H))$.  Then $p \sim
	q$ if and only if the first and the second component
	projections, from the connected component of $(p,q)$ in
	$\cc(\Gamma(H)) \times_A \cc(\Gamma(H))$ to $\cc(\Gamma(H))$
	are both covers.
    \end{itemize}	
\end{proposition}

\proof
The first statement follows directly from Remark~\ref{fact: Lp Lq},
which shows that if $p \sim q$ and there is a path $(p,q) \mapright u
(p',q')$, then $p' \sim q'$.

Let us now assume that $p \sim q$ and let us show that the first
component projection is a cover from the connected component of
$(p,q)$ onto $\cc(\Gamma(H))$.  Let $(r,s)$ be a vertex in that
connected component: then there exists $u\in F$ such that $p \mapright
u r$ and $q \mapright u s$.  Let $r\mapright a r'$ ($a\in \tilde A$)
be an edge in $\Gamma(H)$.  Then $ua \in L_p(H)$, so $ua \in L_q(H)$,
and hence (since $\Gamma(H) \times_A \Gamma(H)$ is deterministic),
there exists an $a$-labeled path $s \mapright a s'$.  Therefore there
exists an $a$-labeled path $(r,s) \mapright a (r',s')$.  Thus the
first component projection is a cover.  The proof concerning the
second component projection is identical.

Conversely, suppose that the first and the second component
projections, from the connected component of $(p,q)$ in
$\cc(\Gamma(H)) \times_A \cc(\Gamma(H))$ to $\cc(\Gamma(H))$ are
covers, and let $u \in L_p(H)$.  Then $\cc(\Gamma(H))$ has a path $p
\mapright u r$.  It is an elementary property of covers that this path
can be lifted to a path in $\cc(\Gamma(H)) \times_A \cc(\Gamma(H))$,
of the form $(p,q) \mapright u (r,s)$.  The second component
projection of that path yields a path $q \mapright u s$ in
$\cc(\Gamma(H))$, and hence $u \in L_q(H)$.
\eop

Let $f(n)$ be the maximal number of finite-index extensions of a
subgroup $H\lefg F$ such that $\Gamma(H)$ has at most $n$ vertices.
By Proposition~\ref{LpLq covers}, every pair $(p,q)$ such that $p\sim
q$ is in the connected component of a pair of the form $(1,r)$ for
some $r > 1$.  Moreover, this connected component has elements of the
form $(i,j)$ for all $1\le i \le n$, so the graph resulting from the
identification of $1$ and $r$ (or from $p$ and $q$) has at most $n/2$
vertices.  Thus $f(1) = 1$ and $f(n) \le n\,f(\lfloor n/2\rfloor)$ for
all $n \ge 2$.  It follows that $f(n) \le n^{\frac12(1+\log_2 n)}$.

\begin{proposition}
    Let $H \lefg F$. If $\cc(\Gamma(H))$ has $n$ vertices, then $H$ 
    has at most $n^{\frac12(1+\log_2 n)}$ finite-index extensions.
\end{proposition}

\begin{example}\label{ex hypercube}
    By means of lower bound, we consider the following example.  Let
    $e_1,\ldots,e_k$ be the canonical basis of the vector space
    $\Z_2^k$, let $\phi$ be the morphism from the free group $F$ over
    $A = \{a_1,\ldots,a_k\}$ into the additive group $\Z_2^k$, mapping
    $a_i$ to $e_i$, and let $H = \ker\phi$.  Then $H$ is normal and
    finite-index, so all its extensions have finite index and they are
    in bijection with the set of quotients of $\Z_2^k$, hence with the
    set of subgroups of $\Z^k_2$, or equivalently with the set of
    subspaces of $\Z_2^k$.
    
    Let $\ell_{d,k}$ be the number of linearly independent $d$-tuples
    in $\Z_2^k$ ($d\ge 1$).  Then $\ell_{1,k} = 2^k-1$.  If $d\ge 2$,
    a $d$-tuple $(x_1,\ldots,x_d)$ is linearly independent if and only
    if $(x_1,\ldots,x_{d-1})$ is linearly independent and $x_d$ does
    not belong to the subspace generated by $x_1,\ldots,x_{d-1}$, so
    that $\ell_{d,k} = \ell_{d-1,k}(2^k - 2^{d-1})$.  Now the set of
    cardinality $d$ linearly independent subsets has $m_{d,k} =
    \ell_{d,k}/d!$ elements, and the number of dimension $d$ subspaces
    of $\Z_2^k$ is
    $$s_{d,k} = \frac{m_{d,k}}{m_{d,d}} =
    \frac{\ell_{d,k}}{\ell_{d,d}} = \frac{(2^k-1)(2^k-2)(2^k-4) \cdots
    (2^k-2^{d-1})}{(2^d-1)(2^d-2)(2^d-4) \cdots (2^d-2^{d-1})}.$$
    Finally, the number of subspaces of $\Z_2^k$ is equal to
    $\sum_{d=0}^k s_{d,k}$, with $s_{0,k} = 1$.
    
    We observe that for each $0 \le i < d < k$,
    $\frac{2^k-2^i}{2^d-2^i} > 2^{k-d}$, so that $s_{d,k} >
    2^{(k-d)d}$.  By considering $d = \lfloor\frac k2\rfloor$, we find
    that $\sum s_{d,k} > 2^{k^2/4}$.
    
    Finally, we note that $\Gamma(H)$ is the Cayley graph of $\Z_2^k$
    with respect to the basis $e_1,\ldots,e_d$ (a graph known as the
    dimension $k$ hypercube), so that $\Gamma(H)$ has $n = 2^k$
    vertices.  As a result, $H$ has more than $n^{\frac14\log_2n}$
    finite-index extensions.
\end{example}


\subsection{The lattice of finite-index subgroups of $G$}\label{sec: 
FIS}

Let us call \textit{\fim-maximal} a subgroup $G \lefg F$ which has no
proper finite-index extension, that is (in view of
Theorem~\ref{maximal fi extension}), such that $G = \Comm_F(G)$.  If
$G$ is \fim-maximal, let $\FIS(G)$ be the set of all finite-index
subgroups of $G$, that is, the set of subgroups $H\lefg F$ such that
$H_\fim = G$.  Note that distinct \fim-maximal subgroups yield
disjoint lattices of finite-index subgroups.

\begin{remark}
    Let $G \lefg F$ be non-trivial and \fim-maximal.  Then $\FIS(G)$
    forms a convex sublattice of the lattice of subgroups of $F$, with
    greatest element $G$.  This sublattice is always infinite (see
    Remark~\ref{remark index}) and without a least element.
\end{remark}

Lemma~\ref{preserve Lr} provides us with an invariant for every
sublattice of the form $\FIS(G)$ (with $G$ \fim-maximal).

\begin{proposition}\label{prop invariant}
    Let $H, K \lefg F$.  Then $H_\fim =
    K_\fim$ if and only if $\tt_H(1) = \tt_K(1)$ and $\widetilde
    L_{\tau(1)}(H) = \widetilde L_{\tau(1)}(K)$, if and only if
    $\tt_H(1) = \tt_K(1)$ and $L_{\tau(1)}(H) = L_{\tau(1)}(K)$.
\end{proposition}    

\proof
If $H_\fim = K_\fim$, then $H,K \lefi H_\fim$, and
Proposition~\ref{prop: fi cover} and Lemma~\ref{preserve Lr} show
directly that $\tt_H(1) = \tt_{H_\fim}(1) = \tt_K(1)$, $\widetilde
L_{\tau(1)}(H) = \widetilde L_{\tau(1)}(H_\fim) = \widetilde
L_{\tau(1)}(K)$ and $L_{\tau(1)}(H) = L_{\tau(1)}(H_\fim) =
L_{\tau(1)}(K)$.

We now prove the converse.  More precisely, we show that if $H$ and
$K$ are \fim-maximal, $\tt_H(1) = \tt_K(1)$ and $L_{\tau(1)}(H) =
L_{\tau(1)}(K)$, then $H = K$.  First we note that
$\calA(H^{\tt_H(1)}) = (\cc(\Gamma(H)),\tau(1))$ and in particular,
$H^{\tt_H(1)}$ is cyclically reduced and \fim-maximal.  Thus, it
suffices to prove the expected result (namely, that $H = K$) under the
hypothesis that $H$ and $K$ are cyclically reduced.

By Remark~\ref{fact: Lp Lq}, the set of all $L_p(H)$ ($p \in \Gamma(H)$)
coincides with the set $\{\{\red(\bar ux) \mid x \in L_1(H)\} \mid u
\in L_1(H)\}$.  In addition, since $H$ is \fim-maximal, if $u,v\in
L_1(H)$ with $1 \mapright u p$ and $1 \mapright u q$ and $\{\red(\bar
ux) \mid x \in L_1(H)\} = \{\red(\bar vx) \mid x \in L_1(H)\}$, then
$p = q$.  It also follows from the same fact that, again if $H$ is
\fim-maximal, there is an edge $(p,a,q)$ in $\Gamma(H)$ if and only if
$L_q(H) = \{\red(\bar ax) \mid x\in L_p(H)\}$.  Thus, the cyclically
reduced \fim-maximal subgroup $H$ is entirely determined by the set
$L_1(H)$.
This concludes the proof.
\eop

The pairs $(t, L)$ that are equal to $(\tt(1),\widetilde
L_{\tau(1)}(G))$ for some subgroup $G \lefg F$ are characterized as
follows.  Recall that an \textit{i-automaton} \cite{SilvaFIM} is a
deterministic automaton $(Q,i,E,T)$ over the alphabet $\tilde A$ such
that, $(p,a,q) \in E$ if and only if $(q,\bar a,p) \in E$ for all
vertices $p,q$ and $a\in A$.  The automata $\calB_p$ discussed in
Section~\ref{sec: computing} are i-automata.

\begin{proposition}
    Let $t \in F$ and let $L \subseteq \tilde A^*$ be a rational
    language.  The following conditions are equivalent.
    \begin{itemize}
	\item[(1)] There exists a subgroup $H \lefg F$ such that
	$t = \tt(1)$ and $L = \widetilde L_{\tau(1)}(H)$.

	\item[(2)] There exists a \fim-maximal subgroup $H \lefg F$ such that
	$t = \tt(1)$ and $L = \widetilde L_{\tau(1)}(H)$.
	
	\item[(3)] $L$ is accepted by an i-automaton with all states
	accepting and such that, for each state $p$, there exist
	transitions $(p,a,q)$ and $(p,b,r)$ for at least two distinct
	letters $a,b \in \tilde A$.  In addition, if $t$ is not the
	empty word, then $t = t'a$ for a letter $a\in \tilde A$ such
	that $\bar a \not\in L$.
	
	\item[(3')] The minimal automaton of $L$ is an i-automaton
	with all states accepting and such that, for each state $p$,
	there exist transitions $(p,a,q)$ and $(p,b,r)$ for at least
	two distinct letters $a,b \in \tilde A$.  In addition, if $t$
	is not the empty word, then $t = t'a$ for a letter $a\in
	\tilde A$ such that $\bar a \not\in L$.
	
	\item[(4)] $t$ and $L$ satisfy the following conditions:
	\begin{itemize}
	    \item[(4.1)] $L$ is closed under taking prefixes;
	    
	    \item[(4.2)] if $u,v \in L$, then $u\bar uv \in L$;
	    
	    \item[(4.3)] if $uv\bar vw \in L$, then $uw\in L$;
	    
	    \item[(4.4)] if $ua \in L$ with $a\in \tilde A$, then $uab \in 
	    L$ for some $b\in \tilde A$ such that $b\ne \bar a$.
	    
	    \item[(4.5)] if $t$ is not the empty word, then $t = t'a$
	    for a letter $a\in \tilde A$ such that $\bar a \not\in L$.
	\end{itemize}
    \end{itemize}
\end{proposition}

\proof
Conditions (1) and (2) are equivalent by Proposition~\ref{prop
invariant}.

Proposition~\ref{prop compute Hfi} shows that (1) implies (3'), which
in turn implies (3).  Let us now assume that (3) holds and let
$\Gamma$ be the $A$-labeled graph induced by the states and
transitions of the minimal automaton of $L$.  The extra condition
given shows that every vertex of $\Gamma$ is visited by a loop labeled
by a cyclically reduced word.  It follows that, if $G$ is the subgroup
whose graphical representation is $\calA(G) = (\Gamma,q_0)$, with
$q_0$ the initial state, then $G$ is cyclically reduced, $L =
\widetilde L_{q_0}(G)$.  The condition on the word $t$ shows that
Condition (1) holds with $H = G^t$.

Condition (3) easily implies Condition (4).  Let us now assume that
Condition (4) holds.  By \cite[Theorem 4.1]{SilvaFIM}, Properties
(4.2) and (4.3) show that $L$ is accepted by an i-automaton.  Property
(4.1) shows that all states of that automaton are final, and Property
(4.4) shows that, for each state $p$, there exist transitions
$(p,a,q)$ and $(p,b,r)$ for at least two distinct letters $a,b \in
\tilde A$. Thus Condition (4) implies Condition (3), which concludes 
the proof.
\eop

\section{Malnormal closure}\label{sec: malnormal}

A subgroup $H$ of $F$ is \textit{malnormal} if $H^g \cap H = 1$ for
each $g\not\in H$.  Malnormality was proved decidable in
\cite{Baumslag}, and a simple decision algorithm was given in
\cite{KM}, based on the following characterization
\cite{KM,Jitsukawa}.

\begin{proposition}\label{charact malnormal}
    Let $H \lefg F$.  Then $H$ is malnormal if and only if every
    connected component of $\Gamma(H) \times_A \Gamma(H)$, except for
    the diagonal complement, is a tree, if and only if, for every
    $p\ne q\in \cc(\Gamma(H))$, $L_p(H) \cap L_q(H)$ is finite.
\end{proposition}

This yields directly an $\O(n^2\log n)$ decision algorithm, where $n$ 
is the number of vertices of $\cc(\Gamma(H))$. It also yields the 
following corollary.

\begin{corollary}
    Let $H\lefg F$. If $H$ is malnormal, then $H$ is \fim-maximal, 
    that is, $H$ has no proper finite-index extension.
\end{corollary}

\proof
By Corollary~\ref{charact fi extensions}, if $H$ is not \fim-maximal,
then there exist vertices $p\ne q$ in $\cc(\Gamma(H))$ such that
$L_p(H) = L_q(H)$, and hence such that $L_p(H) \cap L_q(H) = L_p(H)$
is infinite.  In particular, $H$ is not malnormal.
\eop

It is shown in \cite[Prop.  4.5]{MVW} that for every finitely
generated subgroup $H \lefg F$, there exists a least malnormal
extension $H_\mal$, called the \textit{malnormal closure} of $H$, that
$H_\mal$ is finitely generated and effectively computable, and that
the rank of $H_\mal$ is less than or equal to the rank of $H$
\cite[Corol.  4.14]{MVW}.  In fact, \cite{MVW} shows that
$\Gamma(H_\mal)$ is obtained from $\Gamma(H)$ by a series of i-steps.
The algorithm computing $H_\mal$ then consists in computing all the
quotients of $\Gamma(H)$ and verifying, for each of them, whether it
represents a malnormal subgroup.  We now give a much better,
polynomial-time algorithm.

\begin{theorem}\label{good algo malnormal}
    Let $H \lefg F_n$.  The malnormal closure of $H$ is computed by
    repeatedly applying the following construction: in $\Gamma(H)$,
    identify all pairs $(p,q)$ of distinct vertices in
    $\cc(\Gamma(H))$, such that $L_p(H) \cap L_q(H)$ is infinite and
    reduce the resulting graph.
\end{theorem}

\proof
Let $H_0 = H$ and let $H_{i+1}$ be the subgroup of $F$ such that
$\Gamma(H_{i+1})$ is obtained from $\Gamma(H_i)$ by, first identifying
all pairs $p,q$ of vertices of $\cc(\Gamma(H_i))$ such that $L_p(H_i)
\cap L_q(H_i)$ is infinite, and then reducing the resulting graph.
Since $\Gamma(H_{i+1})$ has less vertices than $\Gamma(H_i)$,
this defines a finite sequence of subgroups
$$H = H_0 < H_1 < \ldots < H_k,$$
where $k \ge 0$ and $H_k$ is malnormal.  We want to show that $H_k$ is
the least malnormal extension of $H$.

Let $K$ be any malnormal subgroup of $F$ such that $H \le K$.  By
Takahasi's theorem (see for instance \cite{MVW}), there exists a
subgroup $G$ such that $H \le G \le K$ such that $G$ is obtained from
$H$ by a sequence of i-steps, and $K$ is obtained from $G$ by a
sequence of re-steps (this fact can also be deduced from \cite[Prop.
2.6]{SW08}).

Let $p,q \in \cc(\Gamma(H))$ such that $L_p(H) \cap L_q(H)$ is
infinite.  It is elementary to verify that $L_p(H) \subseteq
L_{\phi^G_H(p)}(G) \subseteq L_{\phi^K_G(\phi^G_H(p))}(K)$.  In
particular, $\phi^K_G(\phi^G_H(p)) \cap \phi^K_G(\phi^G_H(q))$ is
infinite.  Since $K$ is malnormal, it follows from
Proposition~\ref{charact malnormal} that $\phi^K_G(\phi^G_H(p)) =
\phi^K_G(\phi^G_H(q))$.  But $\phi^K_G$ is one-to-one by definition,
so we have $\phi^G_H(p) = \phi^G_H(q)$.  It follows that $G$ is
obtained from $H_1$ by a sequence of i-steps, and in particular, $H_1
\le G \le K$.  Iterating this reasoning, we find that $H_k \le K$,
which concludes the proof.
\eop

\begin{corollary}
    Let $H \le F$.  If $\Gamma(H)$ has $n$ vertices, then one can
    compute $\Gamma(H_\mal)$ in time $\O(n^3\log n)$.
\end{corollary}    

\proof
According to the algorithm given in Proposition~\ref{good algo
malnormal}, we first need to compute the connected components of
$\cc(\Gamma(H)) \times_A \cc(\Gamma(H))$ -- done in time $\O(n^2\log
n)$ --, identifying which are trees and which are not -- done in time
$\O(n^2)$ --, identifying the vertices of $\cc(\Gamma(H))$ involved in
a non-diagonal connected component and reducing the resulting graph to
obtain $\Gamma(H_1)$ -- which is done in time $\O(n\log n)$.  Thus
$\Gamma(H_1)$ is computed from $\Gamma(H_0)$ in time $\O(n^2\log n)$.

This part of the algorithm is iterated $k$ times, to compute
$\Gamma(H_k) = \Gamma(H_\mal)$, and we have $k < n$ since the number
of vertices of the $\Gamma(H_i)$ forms a properly decreasing sequence.
This concludes the proof.
\eop

\subsection*{Acknowledgements}
The authors thank A. Martino and E. Ventura for their help in 
the computation in Example~\ref{ex hypercube}. They also  are 
indebted to the anonymous referee whose careful reading helped 
improve this paper.

    
\end{document}